\documentclass{amsart}
\usepackage{amsmath}
\usepackage{amscd}
\usepackage{epsfig}
\usepackage{amsthm}
\usepackage{graphicx}
\usepackage{amssymb}
\usepackage{latexsym}
\usepackage{epsfig}
\usepackage{epic}

\newtheorem{theorem}{Theorem}[section]

\def\Re{\mathop{\rm Re}\nolimits}
\def\M{\mathop{\bf M}\nolimits}

\def\Todd{\mathop{\rm Todd}\nolimits}
\begin{document}
\title[The Ehrhart function for symbols]{The Ehrhart function for symbols}
\author[V.W.Guillemin]{Victor W. Guillemin}
\address{Department of Mathematics, Mass. Inst. of Tech., Cambridge, MA 02139}
\email{vwg@math.mit.edu}

\author[S.\ Sternberg]{Shlomo Sternberg}
\address{Department of Mathematics, Harvard University,
Cambridge, MA 02138, USA}
\email{shlomo@math.harvard.edu}

\author[J.\ Weitsman]{Jonathan Weitsman}
\address{Department of Mathematics, University of California,
Santa Cruz, CA 95064, USA}
\email{weitsman@math.UCSC.EDU}
\thanks{The work of J.W. is partially supported by National Science Foundation Grant
DMS/0405670.}
\begin{abstract}
We derive an Ehrhart function for symbols from the Euler-MacLaurin 
formula with remainder.
\end{abstract}
\maketitle
\begin{center} In memory of Prof. S.S. Chern, our teacher and mentor.\end{center}
\tableofcontents
\section{Introduction.}
Let $\Delta\subset\mathbb{R}^n$ be a convex polytope whose vertices are in $\mathbb{Z}^n$ 
and such that the origin $0$ is in the interior of $\Delta$.  Consider the expanded polytope
$N\cdot \Delta$.

Ehrhart's theorem \cite{Ehr} asserts that for $N$ a positive integer, the number of lattice points in 
the expanded polytope, i.e.  
$$\#\left((N\cdot \Delta)\cap\mathbb{Z}^n\right),\ \ \ N\in \mathbb{Z}_+$$
is a polynomial in $N$. More generally, suppose that $f$ is a polynomial, and let
\begin{equation}
\label{definep(N,f)}
p(N,f):=\sum_{\ell\in N\cdot\Delta\cap\mathbb{Z}^n}f(\ell).
\end{equation}
Then Ehrhart's theorem asserts that $p(N,f)$ is a polynomial in $N$. 

In the case that $\Delta$ is a simple polytope (meaning that $n$ edges emanate from 
each vertex) Ehrhart's theorem is a consequence of the Euler-MacLaurin formula, \cite{Kh, KP,CS:bulletin,CS:EM,Gu,BV,DR} and one
can be more explicit about  the nature of the polynomial $p(N,f)$. 

Let us explain how this works in the more restrictive case where $\Delta$ is not only 
simple but is regular, meaning that the local cone at each vertex can be transformed 
by an integral unimodular affine transformation into a neighborhood  of the origin in the standard orthant  
$\mathbb{R}_+^n$.\footnote{Of course, this use of the word ``regular" has 
nothing to do with the term denoting Platonic solids. There are other names in the literature for the property we are describing,  such as   ``smooth", ``Delzant" ,``torsionfree", ``unimodular"
etc.. It is unfortunate that the nomenclature for polytopes with this property has not yet been standardized.} In this case we can apply the formula of Khovanskii-Pukhlikov \cite{KP}, which 
reads as follows: The polytope $\Delta$   can be described by a set of inequalities 
$$x\cdot u_i+a_i\geq 0,\ \ \ \ i=1,\dots,m$$
where $m$ is the number of facets of $\Delta$, where the $u_i$ are primitive lattice vectors,
and where the $a_i$ are positive integers. Then for any positive number $t$, 
the expanded polytope $t\cdot\Delta$ is described by the inequalities 
$$x\cdot u_i+ ta_i\geq 0,\ \ \ \ i=1,\dots,m.$$
Let $\Delta_{t,h}, \ \  h=(h_1,\dots, h_m)$ be the polytope defined by
\begin{equation}
\label{defineDeltath}
x\cdot u_i+ ta_i+h_i\geq 0,\ \ \ \ i=1,\dots,m.
\end{equation}
Then the function
\begin{equation}
\label{definetildapf}
\tilde{p}(t,h,f):=\int_{\Delta_{t,h}}f(x)dx
\end{equation}
is a polynomial in $t$ and $h$. 

The  formula of Khovanskii-Pukhlikov (applied to $N\cdot \Delta$) expresses 
 $p(N,f)$ in terms of a differential operator applied to $\tilde{p}(t,h,f)$: Explicitly, consider 
the infinite order constant coefficient differential operator 
\begin{equation}
\label{Toddoperator}
\Todd\left(\frac{\partial}{\partial h}\right)=\sum_\gamma b_\gamma \left(\frac{\partial}{\partial h}\right)^\gamma 
\end{equation}
where $\sum_\gamma b_\gamma x^\gamma$ is the Taylor series expansion at the origin of the 
 Todd function
$$\Todd(x)=\prod_{i=1}^m\frac{x_i}{1-e^{-x_i}}.$$
The  Khovanskii-Pukhlikov formula says that
$$p(N,f)=\left.\Todd\left(\frac{\partial}{\partial h}\right)\tilde{p}(N,h,f)\right|_{h=0}.$$
Note that since $\tilde{p}$ is a polynomial in $h$ the right hand side really involves 
only a finite order differential operator.

For purposes below it will be convenient to write the Khovanskii-Pukhlikov  formula in the form
\begin{equation}
\label{KHOPHU}
p(N,f)-\tilde{p}(N,0,f)=\left.\left(\Todd\left(\frac{\partial}{\partial h}\right)-\hbox{Id}\right)
\tilde{p}(N,h,f)\right|_{h=0}.
\end{equation}

For simple polytopes there is a more general formula due to \cite{CS:EM,Gu,BV}. Our goal in this 
paper is to prove an analogue of (\ref{KHOPHU}) and its generalizations when the 
polynomial $f$ is replaced by a ``symbol",
a term whose definition we recall from the theory of partial
differential equations:  A smooth function $f \in C^\infty(\mathbb{R}^n)$ 
is called a \textbf{symbol of order} ${r}$
if for every $n$-tuple of non-negative integers
$a :=(a_1,\dots,a_n)$,
there exists a constant $C_a$ such that 
$$|\partial_1^{a_1}\dots \partial_n^{a_n} f(x) | \leq C_a (1 + |x|)^{r - |a|}$$
where $|a| = \sum_i a_i$.  In particular, a polynomial of degree
$r$ is a symbol of order $r$.  
Note that if $f$ is a symbol of order $r$ on $\mathbb{R}^n$ then its
derivatives of order $a$ are in $L^1$ if $r < |a| - n$.

 For  simplicity,
we will restrict ourselves in this paper to {\bf polyhomogeneous symbols}, meaning functions 
$f\in C^\infty(\mathbb{R}^n)$ which admit asymptotic expansions of the form:
\begin{equation}
\label{polyhomsymb}f(x)\sim \sum_{-\infty}^r f_\ell(x)
\end{equation}
for $\|x\|>>0$ where the $f_\ell$ 
 are homogeneous symbols  of degree $\ell$ meaning that  each $f_\ell$ is a 
symbol of order $\ell$  with the property that for $\|x\|$ sufficiently large and 
$t\geq 1$ we have $f_\ell(tx)=t^\ell f_\ell(x)$.  The sum in (\ref{polyhomsymb})  is over a discrete sequence 
of numbers  tending  to $-\infty$. ``Asymptotic"  means that for any $j$
\begin{equation}
\label{defofasymptotic} f(x)-\sum_{\ell= j}^rf_\ell(x)=o(\|x\|^j)
\end{equation}
as $\|x\|\rightarrow\infty$. The number $r$ occurring in (\ref{polyhomsymb}) is again called the {\bf order}
 of the asymptotic series and the collection of functions satisfying (\ref{polyhomsymb}) will 
be called polyhomogeneous symbols of order $r$ and denoted by $S^r$. 

For the sake of exposition in this introduction, 
we will continue to discuss the case where $\Delta$ is regular. 

We will show that if $f$ is a  polyhomogeneous symbol then the function $p(N,f)$ given by 
 (\ref{definep(N,f)}) is a polyhomogeneous symbol in $N$ and its asymptotic 
 expansion in powers of $N$ is  given by a formula similar to  (\ref{KHOPHU}) with two key 
differences:
\begin{enumerate}
\item For symbols,  an infinite number of differentiations occur on the right hand side of 
(\ref{KHOPHU}), i.e. the whole Todd operator must be applied. So (\ref{KHOPHU}) must be
 understood as an asymptotic series, not as an equality.
\item The formula (\ref{KHOPHU}) has to be corrected by adding a constant term $C$ 
to the right hand side, a constant which is zero for the case of a polynomial.
\end{enumerate}

More precisely, we will prove:
\begin{theorem}\label{Ehrhartforsymbolsreg}
Let $\Delta$ be a regular polytope whose vertices lie in $\mathbb{Z}^n$ with $0$ in the interior
of $\Delta$. Let $f\in S^r$ and $N\in \mathbb{Z}_+$ and let $$p(N,f):=\sum_{\ell\in N\cdot\Delta\cap\mathbb{Z}^n}f(\ell).$$
Let $\tilde{p}(t,h,f)$ be defined by (\ref{definetildapf})  so that 
$$\tilde{p}(N,0,f)=\int_{N\cdot \Delta}f(x)dx,\ \ \  $$
Then $p(N,f)-\tilde{p}(N,0,f)$ is a symbol in $N$ and has the asymptotic expansion
$$p(N,f)-\tilde{p}(N,0,f)\sim \left.\left(\Todd\left(\frac{\partial}{\partial h}\right)-\hbox{\rm{Id}}\right)
\tilde{p}(N,h,f)\right|_{h=0}+C$$
where $C$ is a constant. 
\end{theorem}
The constant $C$ is of interest in its own right. It can be thought of as a ``regularized" 
version of the difference
\begin{equation}
\label{unrenomrlaizeddiff}
\sum_{\ell\in \mathbb{Z}^n}f(\ell)-\int_{\mathbb{R}^n}f(x)dx.
\end{equation}
Of course there is no reason why either the sum or the integral in (\ref{unrenomrlaizeddiff}) should converge. But we 
can ``regularize" both as follows:  Define the function $\langle x\rangle$ by
$$\langle x\rangle^2:=1+\|x\|^2.$$
For $s\in \mathbb{C}$ let
$$f(x,s):=f(x)\langle x\rangle^s.$$
We will show that
\begin{equation}
\label{defineC(s)}
C(s):=\sum_{\ell\in \mathbb{Z}^n}f(\ell,s)-\int_{\mathbb{R}^n}f(x,s)dx,
\end{equation}
which is holomorphic for Re $s<<0$, has an analytic continuation to the entire complex 
plane and that the missing constant $C$ on the right hand side of (\ref{KHOPHU})  
is exactly $C(0)$. In particular, the constant  $C$  is independent of the particular polytope
in question. 

This result is related to, and inspired by, a result of Friedlander-Guillemin \cite{FG} on 
``Szego regularizations" of determinants of pseudodifferential operators. In their result, 
as in ours, there is a missing constant which also has to be computed by the above process 
of ``zeta regularization".

Our result is somewhat insensitive to the mode of regularization. In fact, it can be generalized as
follows:Define a ``gauged symbol" \cite{Gugau} to be a function $f(x,s)\in C^\infty(\mathbb{R}^n\times \mathbb{C})$ which depends holomorphically on $s$ and for fixed $s$ is a symbol 
of order Re $s+r$. For example, if $f$ is a symbol,  the function $f(x)\langle x\rangle^s$ introduced above 
is such a gauged symbol.  We make a similar definition of  ``gauged  polyhomogenous symbols". We will prove that if $f(x,s)$ is a gauged polyhomogenous symbol with
$$f(x)=f(x,0)$$
then the function given by (\ref{defineC(s)}) with this more general definition of $f(x,s)$ 
again extends holomorphically from Re $s<<0$ to the entire plane and $C=C(0)$. 

The above results will be proved for the more general case of simple integral polytopes in \S 2. The proof is largely based on the Euler-MacLaurin 
formula with remainder as proved in \cite{KSW}  and motivated by an argument of Hardy on 
``Ramanujan regularization" [Hardy].  Ramanujan's  key idea was to use the classical Euler-MacLaurin formula
in one variable to regularize (\ref{unrenomrlaizeddiff}) by providing ``counter terms" 
in passing to infinity in the difference between sum and integral in one dimension. 
The origin of this method goes back to Euler's continuation of the zeta function past the pole
at $z=1$ and his introduction of what is known today as ``Euler's constant". 

We would like to thank Yael Karshon and Michele Vergne for helpful comments. 

\section{An Ehrhart formula for simple polytopes.}
We continue with the  notation of \S1. So $\Delta\subset\mathbb{R}^n$ is a convex polytope whose vertices are in $\mathbb{Z}^n$ 
and such that the origin $0$ is in the interior of $\Delta$.  We  assume in this section that $\Delta$ 
is  simple, which means that for every vertex $p$ there are exactly $n$ edges emanating from 
$p$ so they lie on rays
\begin{equation}
\label{raysfromp}p+t\alpha_{i,p}
\end{equation}
where the $\alpha_{i,p},\ i=1,\dots,n$ form a basis of  $\mathbb{R}^n$. 
\subsection{The Euler-MacLaurin formula for symbols.}
We want to apply the Euler-MacLaurin formula with remainder, \cite{KSW}. In \cite{KSW} 
one dealt with a weighted sum where points in the interior of the polytope are given 
weight $1$, points on the relative interior of a facet are given weight $w(x):=\frac12$, and, more generally, 
points in the relative interior of a face of codimension $k$ are given weights $w(x):=\frac1{2^k}$.  
The weighted sum  $p_{\frac12}(N,f)$ is  then defined as
$$p_{\frac12}(N,f):=\sum_{\ell\in(N\cdot\Delta)\cap\mathbb{Z}^n}w(\ell)f(\ell).$$
 Theorem 3 of \cite{KSW} gives an Euler-MacLaurin formula with remainder for weighted sums of
symbols. 

More generally, \cite{AW} consider the Euler-MacLaurin formula with remainder for 
more general weightings including the unweighted sum we considered in \S1. 
We refer to equations (28) and (29) in \cite{AW} for the definition of a general weighting, $w$,
 and we will denote the corresponding weighted sum here by $p_w(N,f)$.  
They 
stated their  formula with remainder for smooth functions of compact support, but the 
passage from the case of  smooth functions of compact support to that of symbols is 
exactly the same as in \cite{KSW}.

To formulate  the Euler-MacLaurin formula with remainder we fix a vector $\xi\in \mathbb{R}^n$   
such that
\begin{equation}
\label{xidotanot0}
\alpha_{i,p}\cdot\xi\neq 0\ \ \ \ \forall\ p\ \hbox{and }i.
\end{equation}
We then define
\begin{eqnarray}
\alpha_{i,p}^\sharp&:=&\left\{\begin{array}{rcl}\alpha_{i,p}&\hbox{if}&\alpha_{i,p}\cdot\xi>0\\
						-\alpha_{i,p}&\hbox{if}&\alpha_{i,p}\cdot\xi<0
						\end{array}\right.\label{definealphasharp}\\
						\nonumber\\
(-1)^p&:=&						
\prod_{i=1}^n\frac{\alpha_{i,p}^\sharp\cdot\xi}{\alpha_{i,p}\cdot\xi}\ \ \label{define-1top}\\
\nonumber\\
\hbox{and}\nonumber\\
C_{p,t}&:=&\left\{tp+
\sum_{i=1}^nt_i\alpha_{i,p}^\sharp,\ \ \ t_i\geq0\right\}.
\end{eqnarray}

There is a certain infinite order differential operator $\M$ (depending on the weighting)  in the variables $h_1,\dots,h_m$ 
with constant term $1$ whose truncation of order $k$ is given by the sum in equation (89) 
of \cite{KSW} (for weight $\frac12$) and the sum in equation (56) in \cite{AW}  (for general weights) such that for any
symbol $f$ of order $r$ and $k>n+r$ 
 \begin{equation}
\label{EMremsym}
p_w(N,f)-\tilde{p}(N,0,f)=\left.\left((\M^{[k]})\left(\frac{\partial}{\partial h}\right)-\hbox{Id}\right)\tilde{p}(N,h,f)\right|_{h=0}
+R^k(f,N)
\end{equation}
where $\M^{[k]}$ denotes the truncation of $\M$ at order $k$  and 
\begin{equation}
\label{formofremainder}
R^k(f,N) =\sum_p(-1)^p\int_{C_{p,N}}\left(\sum_{|\gamma|=k}^{|\gamma|=nk}\phi_{\gamma,k}^pD^\gamma f
\right)dx\end{equation}
where the $\phi_{\gamma,k}^p$ are bounded piecewise smooth periodic functions. 
For an  explicit expression for these functions see \cite{KSW} or \cite{AW}. We will not need this here.
If the polytope is regular, and we use unweighted sums, the operator $\M$ is exactly the 
operator $\Todd$ of \S1.  

Notice that  $\M^{[k]}-1$ has no constant term, so $\left.\left((\M^{[k]})\left(\frac{\partial}{\partial h}\right)-\hbox{Id}\right)\tilde{p}(N,h,f)\right|_{h=0} $  involves integration of derivatives of $f$  over faces of the polytope. These faces  
are moving out to infinity as $N\rightarrow \infty$.
A derivative of a homogeneous summand in   (\ref{polyhomsymb}) is itself a homogenous function.
So as $N\rightarrow \infty$, the integral of this derivative over a face of the polytope is a homogenous function in $N$ whose degree depends on the degree of this derivative and the dimension of the face. 
So the  contributions of the polyhomogenous 
terms with sufficiently negative degree in (\ref{polyhomsymb}) will be homogeneous  terms of  high negative degree
in $N$ in (\ref{EMremsym}). By the same token, each homogenous summand in (\ref{polyhomsymb}) will yield a finite number of terms to each order in (\ref{EMremsym}). 
So the left hand side of (\ref{EMremsym}) will be  polyhomogenous symbol in $N$. 

\subsection{The key idea.}
We want to investigate the behavior of the remainder (\ref{formofremainder}) as $N\rightarrow\infty$.

Since the origin $0$ is in the interior of the polytope, it is in the interior of the cone 
generated by the $\alpha_{i,p}$ for any vertex  $p$. This means that when we use 
the $\alpha_{i,p}$ as a basis of $\mathbb{R}^n$ (based at the origin), the point  $p$ has
strictly negative coefficients relative to this basis. When we flip the basis from 
$\alpha_{i,p}$ to $\alpha_{i,p}^\sharp$, the coefficients of those edges which are actually flipped become positive. 

Condition (\ref{xidotanot0}) implies that the function $x\mapsto \xi\cdot x$ has a unique 
minimum on the polytope,  and that this minimum is achieved at a single vertex, $q$. At this vertex, 
no edges are flipped, and at any other vertex $p \neq q$ at least one edge is flipped. So for
$p\neq q$ we have 

\begin{equation}
\label{pintermsofalphasharp}
p=\sum_{i=1}^na_{i,p}\alpha_{i,p}^\sharp,
\end{equation}
where at least one $a_{i,p}> 0$. 
Then the cone $C_{p,t}$ is contained in the half space 
$$\left\{x=\sum_ix_i\alpha_{i,p}^\sharp,\ \ \ x_i\geq a_{i,p}t
\right\}
$$
and so the $p$-th summand in (\ref{formofremainder}) tends to zero as $N\rightarrow \infty$ 
for $p\neq q$.

At the vertex $q$ we have we have $q=\sum_i a_{i,q}\alpha_{i,q}=\sum_i a_{i,q}\alpha_{i,q}^\sharp$ with all the $a_{i,q}<0$ and so the cone $C_{q,t}$ tends to 
the entire space $\mathbb{R}^n$ as $t\rightarrow \infty$. Thus   (\ref{formofremainder}) tends
to 
\begin{equation}
\label{limitofpthterm}
C_k=\int_{\mathbb{R}^n}\left(\sum_{|\gamma|=k}^{|\gamma|=nk}\phi_{\gamma,k}^q D^\gamma f\right)dx.
\end{equation}

 It  follows from (\ref{EMremsym})  that  this limiting value $C_k$ 
is independent of $k$ for $k$ sufficiently large. Let us call this common limit $C=C(f)$.
It also  follows from (\ref{EMremsym}) that $C(f)$ is independent of the choice
of the polarizing vector $\xi$. 

We shall interpret this limiting value $C$ using regularization in the next section, and 
we will find that  $C$ is also independent of the particular polytope we are expanding. 

If $f$ is a polynomial, so that we  choose $k$ to be greater than the degree of $f$, we see from 
(\ref{limitofpthterm}) that $C=0$,  as it must be from the classical  Ehrhart theorem.
\subsection{Application to polyhomogeneous symbols.}
Now suppose that $f$ is a polyhomogenous symbol.
We can apply the above to each summand in the asymptotic series (\ref{polyhomsymb}). 
But notice that if we choose $j$ sufficiently negative, the function 
$$g(x)=g_j(x)=f(x)-\sum_{\ell=j}^rf_\ell(x)$$
occurring on the left hand side of (\ref{defofasymptotic}) will have the property that
both
$$\sum_{\ell\in \mathbb{Z}^n}g(\ell)\ \ \ \ \hbox{and }\int_{\mathbb{R}^n}g(x)dx$$
are absolutely convergent. Furthermore, given any negative number $m$ we can arrange, 
by choosing $j$ sufficiently negative, that
$$p_w(N,g)-\sum_{\ell\in \mathbb{Z}^n}g(\ell)=o(N^m)$$
and 
$$ \tilde{p}(N,0,g)-\int_{\mathbb{R}^n}g(x)dx=o(N^m)$$
so
$$\left[p_w(N,g)-\tilde{p}(N,0,g)\right]-\left[\sum_{\ell\in \mathbb{Z}^n}g(\ell)-
\int_{\mathbb{R}^n}g(x)dx\right]=o(N^m).$$

So if we define 
\begin{equation}
\label{C(f)}
C(f)=\sum_{\ell=j}^rC(f_\ell) + \left[\sum_{\ell\in \mathbb{Z}^n}g_j(\ell)-
\int_{\mathbb{R}^n}g_j(x)dx\right]
\end{equation}
for $j$ sufficiently negative, then $C(f)$ is independent of the choice of $j$.  Furthermore, we see that if $f$ is a polyhomogeneous symbol, we get an asymptotic expansion of the form
$$p_w(N,f)-\tilde{p}(N,0,f) \sim\sum_{\ell=-\infty}^r  \left.\left(\M\left(\frac{\partial}{\partial h}\right)-\hbox{\rm{Id}}\right)
\tilde{p}(N,h,f_\ell)\right|_{h=0}+C(f),$$
where each level in the asymptotic expansion in $N$ involves only finitely many $f_\ell$. 
By abuse of language, we shall denote this equation as
\begin{equation}
\label{abuse}p_w(N,f)-\tilde{p}(N,0,f)  \sim \left.\left(\M\left(\frac{\partial}{\partial h}\right)-\hbox{\rm{Id}}\right)
\tilde{p}(N,h,f)\right|_{h=0}+C(f).
\end{equation}

\subsection{Regularization.}
Suppose that we replace  $f$ by a gauged polyhomogenous symbol $f(x,s)$ with $f(x,0)\in S^r$. Then the remainder
 term (\ref{formofremainder}) applied to a summand in the asymptotic expansion of $f_s=f(\cdot,s)$ is well defined if Re $s<-r-n+k$. Moreover, if $p\neq q$ so that $a_{i,p}>0$ for some $i$ 
 the $p$-th summand on the right of (\ref{formofremainder}) is of order $O(N^{\Re s+r+n-k})$.

At the
 unique vertex $q$ where no edges are flipped
 $q$-th summand of  (\ref{formofremainder})  differs from the integral 
\begin{equation}
\label{limitofpthtermgauged}
\int_{\mathbb{R}^n}\left(\sum_{|\gamma|=k}^{|\gamma|=nk}\phi_{\gamma,k}^qD_x^\gamma f_\ell (x,s)\right)dx
\end{equation}
by a term of order $O(N^{\Re s+\ell+n-k})$. Thus the gauged version of (\ref{abuse}) is
\begin{equation}
\label{EMremsymgauged}
p_w(N,f_s)-\tilde{p}(N,0,f_s)=\left.\left(\M^{[k]}-\hbox{Id}\right)\left(\frac{\partial}{\partial h}\right)\tilde{p}(N,h,f_s)\right|_{h=0}
+ C_k(s)
\end{equation}
$$\ \ \ \ \ \ +O(N^{\Re s+r+n-k})$$
for
$$\Re s<-n-r+k$$
where $f_s=f(\cdot,s)$ where $C_k(s)$ is(\ref{abuse})  with $f$ replaced by $f_s$, and we 
have computed $C(s)$ by going out to level $k$ in the Euler-MacLaurin expansion. 

All the terms on the right of (\ref{EMremsymgauged})
are holomorphic on the half-plane 
$\Re s<-n-r+k$. 

Letting $k\rightarrow \infty$ we conclude that on this half-plane we have, in the notation of (\ref{abuse}),  
\begin{equation}
\label{EMremsymgaugedlim}
p_w(N,f_s)-\tilde{p}(N,0,f_s)\sim \left.\left(\M\left(\frac{\partial}{\partial h}\right)-\hbox{Id}\right)
\tilde{p}(N,h,f_s)\right|_{h=0} +C(s)
\end{equation}
where
\begin{equation}
\label{limitofCkisC}
C(s)-C_k(s)=O(N^{\Re s+r+n-k}).
\end{equation}
Since the $C_k(s)$ are holomorphic of the half-plane $\Re s<-n-r+k$, it follows that 
$C(s)$ is holomorphic on the whole plane.

Moreover, in the asymptotic series on the right of (\ref{EMremsymgaugedlim}) all the terms are
of order at most $\Re s+r+n$. Hence for $\Re s<-r-n$ these terms tend to zero and  we get
\begin{eqnarray*}
C(s)&=&\lim_{N\rightarrow \infty}\left(p(N,f_s)-\tilde{p}(N,0,f_s)\right)\\
&=&\sum_{\ell\in \mathbb{Z}^n}f(\ell,s)-\int_{\mathbb{R}^n}f(x,s)dx,
\end{eqnarray*}
and both the sum and the integral converge absolutely. So if we set $s=0$ we obtain
\begin{equation}
\label{EMremsymgaugedlimfinal}
p_w(N,f)-\tilde{p}(N,0,f)\sim \left.\left(\M\left(\frac{\partial}{\partial h}\right)-\hbox{Id}\right)
\tilde{p}(N,h,f)\right|_{h=0} +C
\end{equation}
where $f(x)=f(x,0)$ and $C=C(0)$. So we can think of $C$ as a ``regularization" of 
 (\ref{unrenomrlaizeddiff}).
To summarize: We have proved
\begin{theorem}\label{Ehrhartforsymbolsrsimple}
Let $\Delta$ be a simple polytope whose vertices lie in $\mathbb{Z}^n$ with $0$ in the interior
of $\Delta$. Let $f\in S^r$ and $N\in \mathbb{Z}_+$ and let $$p_w(N,f):=\sum_{\ell\in N\cdot\Delta\cap\mathbb{Z}^n}w(\ell)f(\ell).$$
Let $\tilde{p}(t,h,f)$ be defined by (\ref{definetildapf})  so that 
$$\tilde{p}(N,0,f)=\int_{N\cdot \Delta}f(x)dx,\ \ \  $$
Then $p(N,f)-\tilde{p}(N,0,f)$ is a symbol in $N$ and has the asymptotic expansion
$$p_w(N,f)-\tilde{p}(N,0,f)\sim \left.\left(\M\left(\frac{\partial}{\partial h}\right)-\hbox{\rm{Id}}\right)
\tilde{p}(N,h,f)\right|_{h=0}+C$$
where $C$ is a constant. Furthermore, if $f(x,s)$ is a gauged polyhomogenous symbol with $f(x,0)=f(x)$ 
then $C=C(0)$ where $C(s)$ is the entire function given by (\ref{EMremsymgaugedlim}) 
and (\ref{limitofCkisC}). For $\Re s<-r-n$
$$C(s)=\sum_{\ell\in \mathbb{Z}^n}f(\ell,s)-\int_{\mathbb{R}^n}f(x,s)dx.$$
Hence $C(s)$ and in particular $C=C(0)$ is independent of the polytope.
\end{theorem}

\noindent{\bf Remarks. 1.} In the course of the discussion we have proved a similar theorem with 
polyhomogeneous symbols replaced by symbols and gauged polyhomogenous symbols
replaced by gauged symbols.

{\bf 2.} Since the initial posting of this paper we have received the interesting paper 
\cite{MP}.

\end{document}